\documentclass{article}
\usepackage{graphicx,graphics}
\usepackage{amsmath,amsthm,amssymb}

\usepackage[all,cmtip]{xy}

\newtheorem{theorem}{Theorem}[section]

\theoremstyle{remark}

\begin{document}

\title{A nerve lemma for gluing together incoherent discrete Morse functions}

\author{Alexander Engstr\"om  \\ Aalto University \\ Helsinki, Finland \\
{\tt alexander.engstrom@aalto.fi}}
\date\today

\maketitle

\begin{abstract}
Two of the most useful tools in topological combinatorics are
the nerve lemma and discrete Morse theory. In this note
we introduce a theorem that interpolates between them
and allows decompositions of complexes into non-contractible
pieces as long as discrete Morse theory ensures that
they behave well enough. The proof is based on diagrams of spaces, 
but that theory is not needed for the formulation or
application of the theorem.
\end{abstract}

\section{Introduction}

Two workhorses of topological combinatorics are the nerve lemma
and discrete Morse theory. In this note we present a formula that
interpolates between them, and allows you to glue together incoherent
discrete Morse functions.

To state the main theorem we need one small new piece of notation.
When recursively building the Morse complex by gluing in critical cells,
the first critical cell considered is the \emph{initial critical cell}. When
there are several critical vertices there is a choice to be made, but the following
theorem doesn't depend on that.

\begin{theorem}\label{theTheorem}
Let $X$ be a  regular CW complex with subcomplexes $X_1, \ldots, X_n$ 
satisfying $X=\cup_{i=1}^n X_i,$ and let $P=\{ \cap_{i\in I} X_i \mid  I \subset \{1,\ldots,n \}\}$
be  a poset with $Y \leq Y'$ whenever $Y \supseteq Y'.$ Fix a Morse matching 
 on each $Y \in P.$ 
If $\dim \sigma_Y > \dim \sigma_Z$ for all comparable subspaces $Y\supset Z$ in $P$
with non-initial critical cells $\sigma_Y \in Y$Êand $\sigma_Z  \in Z,$ then 
\[
X \simeq \bigvee_{Y \in P} Y \ast \Delta(P_{<Y}).
\]
\end{theorem}

Theorem~\ref{theTheorem} interpolate between the
nerve lemma and discrete Morse theory: If all involved spaces are
collapsible, then it's essentially the Nerve lemma; and if there is only one subcomplex,
then it's ordinary discrete Morse theory. In the next section we give an example of how to apply
the formula, and in the last section we prove it.

\section{An example}

The vertices of the
$4\times4$ chessboard complex
\raisebox{-0.3ex}{\makebox[2ex][l]{\includegraphics[height=2ex]{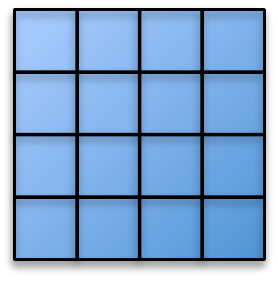}}}
are the $4\times 4$ different places one could put a rook on the board, and
the $4!$ facets correspond to the ways one can place four rooks without having
them attacking each other. In general, the vertices of $n \times n$ chessboard complexes
also correspond to bijections between $n$ elements, and there is a huge
theory on topological, combinatorial and representation theoretic results.
See for example \cite{shareshianWachs} for a survey.

\begin{figure}
 \begin{center}
  \includegraphics[width=0.7\textwidth]{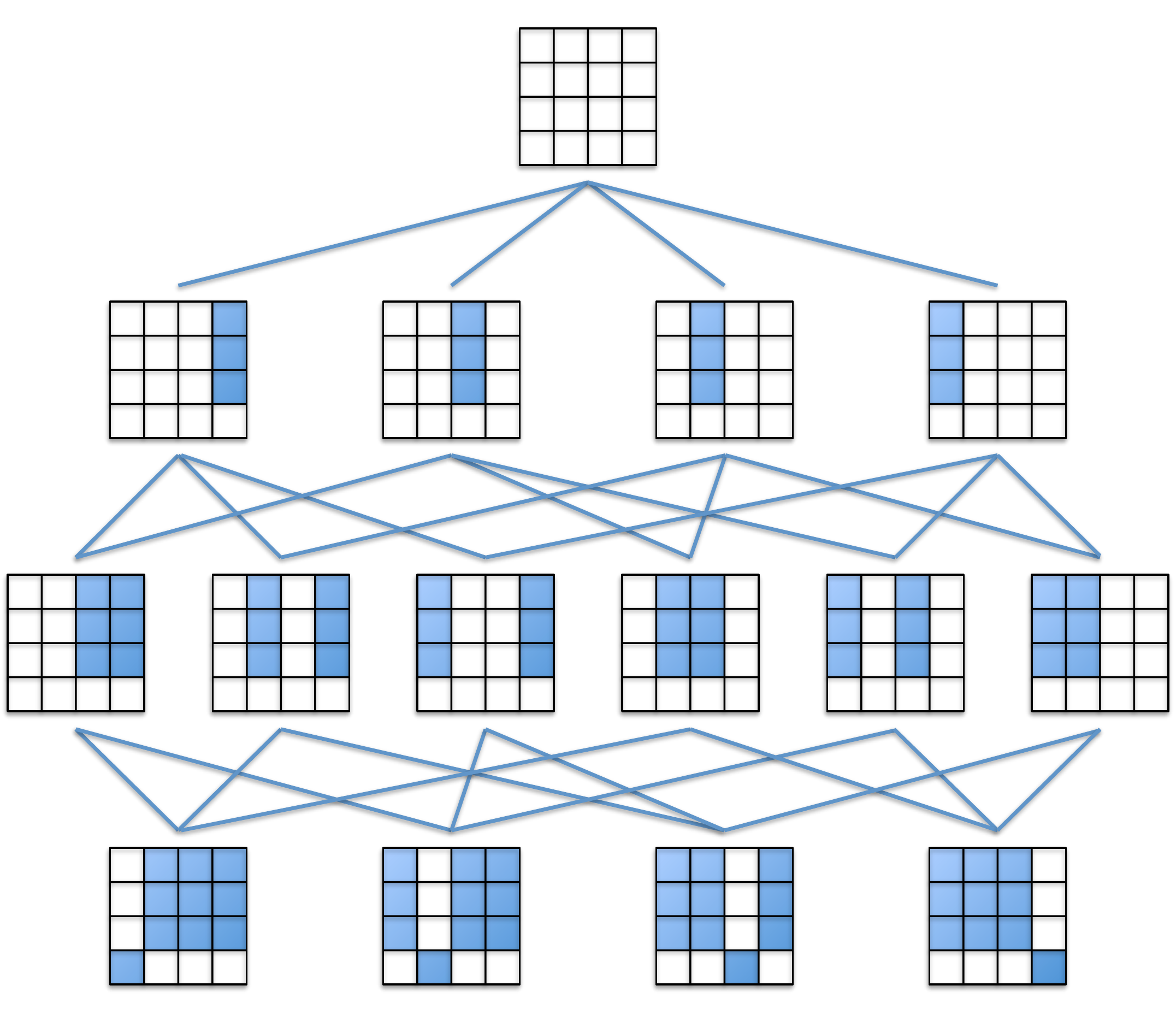}
 \caption{The intersection poset for calculating the homotopy type of a $4\times 4$ chessboard complex.}\label{fig:cb}
 \end{center} 
 \end{figure}

Now we define four subcomplexes of \raisebox{-0.3ex}{\makebox[2ex][l]{\includegraphics[height=2ex]{44_all.pdf}}}
as the vertex induced subcomplexes defined by the shadowed positions in the
bottom row of the poset in Figure~\ref{fig:cb}.
The topological situations for different intersections are as follows:
\begin{itemize}
\item[(Y=\raisebox{-0.95ex}{\makebox[2.4ex][l]{\includegraphics[height=3ex]{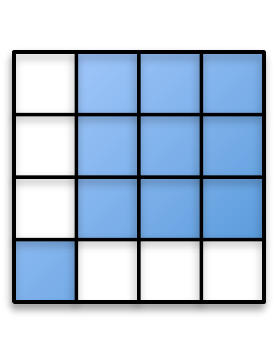}}}) ] The complex $Y$ is a star of a vertex, and a cone with that vertex as apex. Every cone
has a complete Morse matching with only one critical cell, the initial one. Thus, $Y \ast \Delta(P_{<Y}) \simeq \cdot \ast \emptyset \simeq \cdot.$
\item[(Y=\raisebox{-0.95ex}{\makebox[2.4ex][l]{\includegraphics[height=3ex]{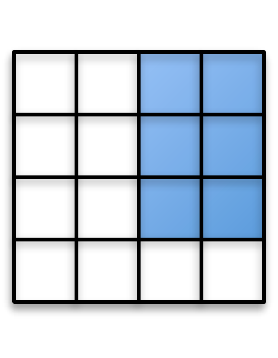}}}) ] The complex $Y$ is an $S^1$ realized as cycle on six edges and six vertices. There is a Morse function with one non-initial critical cell of dimension one, and one initial critical cell of dimension zero. The order complex $\Delta(P_{<Y})$ is two disjoint points, denoted $S^0.$ Thus, $Y \ast \Delta(P_{<Y}) = S^1 \ast S^0 \simeq S^2.$
\item[(Y=\raisebox{-0.95ex}{\makebox[2.4ex][l]{\includegraphics[height=3ex]{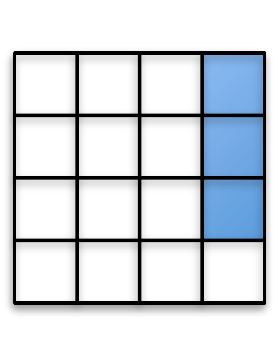}}}) ] The complex $Y$ consists of three disjoint vertices, and every discrete Morse function makes them all critical. The poset $P_{<Y}$ is the face poset of the boundary of a triangle, and $\Delta(P_{<Y})$ is its barycentric subdivision.
$Y \ast \Delta(P_{<Y}) = ( \cdot \cdot \cdot  ) \, \ast S^1  \simeq S^2 \vee S^2.$
\item[(Y=\raisebox{-0.95ex}{\makebox[2.4ex][l]{\includegraphics[height=3ex]{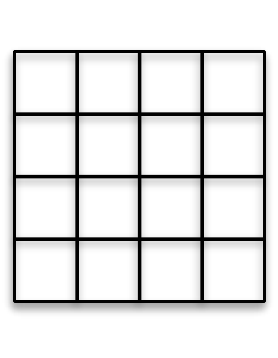}}}) ] The complex $Y$ is empty. The poset $P_{<Y}$ is the face poset of the boundary of a tetrahedron, and $\Delta(P_{<Y})$ is its barycentric subdivision. Thus, $Y \ast \Delta(P_{<Y}) =  \emptyset  \ast S^{2} \simeq S^2.$
\end{itemize}
Compiling the facts, we get that \raisebox{-0.3ex}{\makebox[2ex][l]{\includegraphics[height=2ex]{44_all.pdf}}} is homotopy equivalent to a wedge of 15 two-spheres. This result is well-known, but the computations were sometimes more cumbersome.

\section{Proof of the formula}

Before proving the formula, we point out some references to the methods used in it. As a general
reference for topological combinatorics, Bj\"orner's survey \cite{bjorner} is a good start. The
foundational paper on discrete Morse theory was by Forman \cite{forman}, and for a survey
on contemporary usage, see the book by Jonsson \cite{jonsson}. We build on the explicit
proof of the main theorem of discrete Morse theory provided in Section 3 of \cite{kozlov}.
For diagrams of spaces, the paper by Welker, Ziegler and \v{Z}ivaljevi\'c \cite{welkerZieglerZivaljevic} is recommended.
All of these tools were also used by Engstr\"om to find topological representations of 
matroids by homotopy colimits, and are briefly surveyed 
in \cite{engstrom}. 

\begin{proof}(of Theorem~\ref{theTheorem})
Construct an intersection diagram $\cal{D}$ based on $X=\cup_{i=1}^n X_i.$ This is a diagram
on the poset $P,$ but the language of diagrams of spaces was avoided in the theorem statement.
Intersection diagrams are cofibrant, so $X = \textrm{colim}\,\, \cal{D} \simeq  \textrm{hocolim}\,\, \cal{D}.$

For a space $Y$ with a discrete Morse matching, there are many homotopy equivalences 
$f:Y \rightarrow Y^M$ to the Morse complex. We want to employ a particular type. Let
$y=Y_1 \subset Y_2 \cdots \subset Y_m = Y$ be a filtration such that $y$ is the
initial critical vertex and in each step either two matched cells or a critical cell is added.
This induces a filtration on the Morse complex $y^M=Y_1^M \subset Y_2^M \cdots \subset Y_m^M = Y^M$
and a sequence of homotopy equivalences $f_i:Y_i \rightarrow Y_i^M$ satisfying:
\begin{itemize}
\item[(i)] If $Y_{i}$ extends $Y_{i-1}$ by a critical cell, then $f_i$ extends $f_{i-1}$ by an
homeomorphism from the open cell  $Y_{i} \setminus Y_{i-1}$ to the open cell 
 $Y_{i}^M \setminus Y_{i-1}^M;$
 \item[(ii)] If $Y_{i}$ extends $Y_{i-1}$ by a matching the open cells $\partial \tau \supset \sigma,$ then
 for some homotopy equivalence $\pi: \overline{\tau} \rightarrow \partial \tau \setminus \sigma$ that
 is constant on its image, $f_{i-1}$ extends by setting $f_i(x)=f_{i-1}(\pi(x))$ for 
 $x\in \overline{\tau} \setminus (\partial \tau \setminus \sigma).$
\end{itemize}
A crucial property for the homotopy equivalence $f: Y \rightarrow Y^M$ is that if $x\in Y$ is in the
interior of a $d$-dimensional cell and $f(x)$ is in the interior of a $d^M$-dimensional cell, then
$d\geq d^M.$

For each space $Y$ in $\cal{D}$ a discrete Morse function is given, together with a Morse complex
$Y^M,$ a homotopy equivalence $f_Y:Y \rightarrow Y^M,$ and an initial critical cell $y$ with its image $y^M=f_Y(y).$ 
For any two comparable spaces $Y\supset Z$ in $\cal{D},$ the diagram
\[ \xymatrix{
{}_{\,\,\,}Z \ar@{^{(}->}[d] \ar[r]^{f_Z} & Z^M \ar[d]^{z \mapsto y^M} \\
{\,}Y \ar[r]^{f_Y} &Y^M
} \]
commutes, since the dimension of any non-initial critical cell of $Z$ is less than the dimension
of any non-initial critical cell of $Y.$ Construct a diagram ${\cal D}^M$ on the same poset as $\cal{D}$
but replace each complex by its Morse complex, and each map by a constant map to the image of
the initial critical cell in the corresponding Morse complex. By commutativity of the preceding diagram,
$\textrm{hocolim}\,\, \cal{D} \simeq \textrm{hocolim}\,\, {\cal D}^M.$ Every map in $\cal{D}^M$ is constant,
so $\textrm{hocolim}\,\, {\cal D}^M  \simeq \vee_{Y \in P} Y \ast \Delta(P_{<Y}).$
\end{proof}

\subsubsection*{Acknowledgments.}
This note was prepared for a summer school on discrete Morse theory and commutative algebra,
organized by Bruno Benedetti and the author, July 18 -- August 2, 2012, at Institut Mittag-Leffler.
We thank the institute for being the perfect hosts for our school.

\end{document}